\documentstyle[12pt]{article}
\author{A.A. Balinsky\thanks{ Research supported by
Australian Research Council Grant ( Chief investigator:
\mbox{Prof. G.I. Lehrer )} } \  and
Yu.M. Burman }
\title{Quadratic Poisson brackets and Drinfeld theory for associative
algebras}
\date{}
\begin{document}
\maketitle

%ADD.STY Version 2.0

%New theorem-like environments
% ------------------Numbered-------------------

\newtheorem {theorem}{Theorem}
\newtheorem {lemma}{Lemma}
 \newtheorem {corollary}{Corollary}
%\@addtoreset{corollary}{theorem}
%\@addtoreset{corollary}{lemma}
\newtheorem {statement}{Statement}
\newcounter {rem}
\newenvironment {remark} {\refstepcounter{rem}\par\medskip\noindent{\bf
Remark \arabic{rem}\,\,} }{\par}
\newcounter {exa}
\newenvironment {example} {\refstepcounter{exa}\par\medskip\noindent{\bf
Example \arabic{exa}\,\,} }{\par}
\newcounter {def}
\newenvironment {definition} {\refstepcounter{def}\par\medskip\noindent{\bf
Definition \arabic{def}\,\,}} {\par}

% -----------------Not numbered---------------

\newtheorem {Theorem}{Theorem}
\newtheorem {Lemma}{Lemma}
\newtheorem {Corollary}{Corollary}
\newtheorem {Statement}{Statement}
\newenvironment {Remark} {\par\medskip\noindent{\bf Remark\,\,}} {\par}
\newenvironment {Example} {\par\medskip\noindent{\bf Example\,\,}} {\par}
\renewcommand \theTheorem {}
\renewcommand \theLemma {}
\renewcommand \theCorollary {}
\renewcommand \theStatement {}
\newenvironment {proof} {\par\medskip\noindent{\bf Proof\,\,}}
{\qed\par}
\newenvironment {Proof} {\par\medskip\noindent{\bf Proof\,\,}}
{\par}
\newenvironment {Definition} {\par\medskip\noindent{\bf
Definition\,\,}}{\par}

% ------------------ Apparent symbols -----------------

\def \Real {I\!\!R}
\def \Integer {Z\!\!\!Z}
\def \Complex {\,\prime\mskip-2.5\thinmuskip C}
\def \Rational {\,\prime\mskip-2.5\thinmuskip Q}

% ------------------ Norm and absolute value ------------

\def \lnorm#1\rnorm {\vphantom{#1}\left\|\smash{#1}\right\|}
\def \lmod#1\rmod {\vphantom{#1}\left|\smash{#1}\right|}

% ------------------- More symbols ----------------------

\newcommand \bydef {\stackrel{\mbox{\scriptsize def}}{=}}
\newcommand \qed {\,\rule[-.23ex]{1.6ex}{1.6ex}}
\newcommand \Planck {{h^{\mskip-4.5\thinmuskip{-}}}}

% ------------------ Style conventions ------------------

\newcommand \eps {\varepsilon}
\renewcommand \phi {\varphi}
\renewcommand \rho {\varrho}

\font\eux=eufm10
\font\eus=eufm7
\font\euss=eufm5
\newfam\eufam
\textfont\eufam=\eux
\scriptfont\eufam=\eus
\scriptscriptfont\eufam=\euss
\def\frak{\fam\eufam}

\newcommand \Symm {\mathop{\rm Symm}}
\newcommand \ad {\mathop{\rm ad}\nolimits}
\newcommand \CC {C^\infty}
\newcommand {\SchB}[1]{[[#1,#1]]}
\newcommand {\ii} {{\bf i}}
\newcommand {\jj} {{\bf j}}
\newcommand {\kk} {{\bf k}}

\begin{abstract}
The paper is devoted to the Poisson brackets compatible with multiplication
in associative algebras. These brackets are shown to be quadratic and their
relations with the classical Yang--Baxter equation are revealed. The paper
also contains a description of Poisson Lie structures on Lie groups whose
Lie algebras are adjacent to an associative structure.
\end{abstract}

\section{Introduction}
In his famous work \cite{DR1} V.G.Drinfeld found remarkable relations
between skew-symmetric solutions of classical Yang--Baxter equation
%*
\begin{equation}\label{CYBE}
[r^{12},r^{13}] + [r^{12},r^{23}] + [r^{13},r^{23}] = 0
\end{equation}
%*
where $r$ is an element of the tensor square of some Lie algebra $\frak G$,
and Poisson Lie structures on the corresponding Lie group $G$. More
exactly, fix some basis $e_i$ in $\frak G$, and let $r = r^{ij}e_i \wedge
e_j$ be a solution of (\ref{CYBE}). Let us identify $\frak G$ with the
tangent space of $G$ in the unit, and let $E_i$ and $E'_i$ be left- and
right-invariant vector fields on $G$ extending the basis $e_i$.  The
Poisson Lie tensor (see below exact definitions, as well as explanation of
notation) at the point $x \in G$ is then given by
%*
\begin{equation}\label{DrinPoiss}
\pi^{ij}(x) = r^{ij} (E_i(x) \otimes E_j(x) - E'_i(x) \otimes E'_j(x)).
\end{equation}
%*

Formula (\ref{DrinPoiss}) has many interesting consequences (see
\cite{BAL1} for some examples). But computations involving it necessarily
include direct formulas for $E_i(x)$ and $E'_i(x)$ and may be very
cumbersome. Things improve, however, in case Lie algebra $\frak G$ comes
from an {\em associative} one.

If $A$ is an associative algebra, then one may define a Lie algebra $A_L$
taking $[a,b] = ab - ba$. If $A$ has a unit, then the group $G$ of
invertible elements of algebra $A$ is a Lie group whose Lie algebra is
$A_L$. So, to obtain Poisson Lie structures on $G$, one can start from
Poisson brackets compatible with the algebra $A$. These brackets are
relatively easy to study because algebra $A$ is topologically trivial (a
linear space). It turns out to be that, under some conditions of
smoothness, all these brackets are quadratic. Quadratic Poisson brackets
were studied by many authors (see e.g. a series of works
\cite{KUP3,KUP4,KUP5} by B.Kupershmidt who considered an important case of
a full matrix algebra, and also M.A.Semenov-Tyan-Shanskii's \cite{STS1}).
Quadratic brackets compatible with an associative algebra reveal some new
relations to the classical Yang--Baxter equation which suggest a
possibility to obtain their quantization. This can be done using quadratic
algebras (see \cite{Manin}) and will be a subject of the forthcoming paper.
Quantization of the corresponding Poisson Lie groups may also be simpler
than in general case (see e.g.  \cite{DR}) because such Poisson Lie groups
bear a preferred coordinate system (linear functions on the original
algebra).

The paper is structured as follows. Section \ref{Algebras} deals with
quadratic brackets compatible with an associative algebra structure of $A$.
It is shown that these brackets correspond to {\em differentiations} of the
algebra $\Symm(A \otimes A)$ (symmetric elements of $A \otimes A$) with
values in $A \wedge A$. Jacobi identity means that this differentiation
satisfies some version of classical Yang--Baxter equation. Then, it is
shown in Section \ref{Groups} that restriction of the brackets considered
in Section \ref{Algebras} to the group of invertible elements gives rise to
a Poisson Lie structure. It is also proved that any coboundary Poisson Lie
structure on a simply connected Lie group is a Poisson covering over a
Poisson Lie group whose Poisson bracket is quadratic in some global
coordinate system. Section \ref{Examples} contains examples.

This work was finished during a visit by the first author to the
School of Mathematics and Statistics at the University of Sydney.
A.B. would like to thank Professor G.I. Lehrer for his hospitality.
Authors are grateful to J. Donin, B. Kupershmidt, A. Stolin, and N. Zobin
for many stimulating discussions.

\subsection*{Definitions and notation}

Throughout this paper $A$ stands for a finite-dimensional associative
algebra, $e_i\ (i = 1,\dots,n)$ is its (additive) basis, and $a_{ij}^k$ are
structure constants of the algebra with respect to the basis $e_i$:
%*
\begin{equation}\label{StructConst}
e_i e_j = a_{ij}^k e_k
\end{equation}
%*
(a summation over repeated indices will be always assumed). The symbol
$A_L$ will denote an {\em adjacent} Lie algebra, i.e. Lie algebra on the
same space with a commutator $[a,b] = ab - ba$.

As usual, {\em Poisson bracket} $\{\cdot,\cdot\}$ on a smooth manifold $M$
is understood as a Lie algebra structure on the space of smooth functions
$\CC(M)$ satisfying the Leibnitz identity, i.e. a bilinear operation
$\{\cdot,\cdot\}$ such that
\begin{enumerate}
\item $\{f,g\} = -\{g,f\}$,
\item $\{\{f,g\},h\} + \{\{g,h\},f\} + \{\{h,f\},g\} = 0$ (Jacobi
identity), and
\item $\{fg,h\} = f\{g,h\} + \{f,h\}g$ (Leibnitz identity).
\end{enumerate}
If Jacobi identity is not required one speaks about pre-Poisson brackets.
One can easily see that in local coordinates $x^i$ on the manifold $M$ an
arbitrary pre-Poisson bracket looks like
%*
\begin{equation}\label{PoiTens}
\{f,g\}(x) = \pi^{ij}(x)\partial f(x)/\partial x^i \cdot \partial
g(x)/\partial x^j,
\end{equation}
%*
for some {\em Poisson tensor} $\pi^{ij}(x)$. Bracket (\ref{PoiTens}) is
Poisson if and only if the Poisson tensor satisfies the equation:
%*
\begin{equation}\label{JaPoi}
\pi^{lk}(x) \frac{\partial \pi^{ij}(x)}{\partial x^l} + \pi^{li}(x)
\frac{\partial \pi^{jk}(x)}{\partial x^l} + \pi^{lj}(x) \frac{\partial
\pi^{ki}(x)}{\partial x^l} = 0.
\end{equation}
%*

The set of all Poisson brackets on a given manifold is not a linear space
but just a cone. Indeed, a sum of two Poisson brackets is pre-Poisson, but
it usually does not satisfy Jacobi identity. We call two brackets
$\{\cdot,\cdot\}_1$ and $\{\cdot,\cdot\}_2$ {\em compatible} with one
another if any linear combination $\alpha \{\cdot,\cdot\}_1 + \beta
\{\cdot,\cdot\}_2$ is again a Poisson bracket.

A product $M_1 \times M_2$ of two manifolds equipped with Poisson brackets
$\{\cdot,\cdot\}_1$ and $\{\cdot,\cdot\}_2$, respectively, may be given a
{\em product} Poisson bracket. The latter is a Poisson bracket whose
Poisson tensor $\pi_p^{ij}(x,y)$ in the point $(x,y) \in M_1 \times M_2$ is
%*
\begin{equation}\label{Product}
\pi_p(x,y) =
\left( \begin{array}{cc}
\pi_1^{ij}(x) & 0 \\
0 & \pi_2^{ij}(y)
\end{array}\right)
\end{equation}
%*
where $\pi_1^{ij}$ and $\pi_2^{ij}$ are Poisson tensors of
$\{\cdot,\cdot\}_1$ and $\{\cdot,\cdot\}_2$, respectively. If $p_1 : M_1
\times M_2 \to M_1$ and $p_2 : M_1 \times M_2 \to M_2$ are natural
projections, and $p_1^* : \CC(M_1) \to \CC(M_1 \times M_2)$ and $p_2^* :
\CC(M_2) \to \CC(M_1 \times M_2)$ are corresponding pullbacks, then it is
easy to see that product Poisson bracket may be defined as the only Poisson
bracket on $M_1 \times M_2$ satisfying the conditions:
\begin{enumerate}
\item $\{p_1^*(f),p_1^*(g)\} = \{f,g\}_1$, and $\{p_2^*(f),p_2^*(g)\} =
\{f,g\}_2$,
\item $\{p_1^*(f),p_2^*(g)\} = 0$.
\end{enumerate}

A mapping $F:M_1 \to M_2$ of two manifolds equipped with Poisson brackets
$\{\cdot,\cdot\}_1$ and $\{\cdot,\cdot\}_2$, respectively, is called {\em
Poisson} if a pullback $F^*: \CC(M_2) \to \CC(M_1)$ is a Lie algebra
homomorphism, i.e.
%*
\begin{displaymath}
\{f\circ F,g \circ F\}_1 = \{f,g\}_2 \circ F.
\end{displaymath}
%*

Let $M_1$ be a manifold with a Poisson bracket $\{\cdot,\cdot\}$, and $M'
\subset M$ be its submanifold. Then $M'$ is called {\em Poisson
submanifold} if it carries a Poisson bracket $\{\cdot,\cdot\}'$ such that
the inclusion $\imath : M' \hookrightarrow M$ is a Poisson mapping. It is
easy to see that if the bracket $\{\cdot,\cdot\}'$ exists, then it is
unique. We call it a {\em restriction} of the Poisson bracket
$\{\cdot,\cdot\}$ to the submanifold $M'$.

Let now manifold $M$ be a Lie group, i.e. be equipped with a smooth
multiplication $*: M \times M \to M$. A Poisson bracket on $M$ is called
{\em compatible} with this multiplication if $*$ is a Poisson mapping (with
respect to a product Poisson structure in $M \times M$). In this case $M$
is called a {\em Poisson Lie group}. More explicitly compatibility means
the following: let $f,\,g$ be functions from $\CC(M)$, and $r(x) =
\{f(x),g(x)\}$ be their Poisson bracket. Let also $F,\,G,\,R$ be functions
from $\CC(M \times M)$ given by:
%*
\begin{eqnarray*}
F(x,y) &=& f(x*y), \\
G(x,y) &=& g(x*y), \\
R(x,y) &=& r(x*y).
\end{eqnarray*}
%*
Then there should be
%*
\begin{equation}\label{Compat}
R = \{F,G\}_p
\end{equation}
%*
where $\{\cdot,\cdot\}_p$ means a product Poisson bracket (\ref{Product}).

\section{Quadratic Poisson brackets compatible
  with an algebra structure}\label{Algebras}
Consider now Poisson structures on a linear space $A$. The dual space $A^*$
is naturally embedded into functions algebra $\CC(A)$, as well as all its
symmetric powers $\Symm\bigl((A^*)^{\otimes n}\bigr)$ are. In view of
(\ref{PoiTens}) speaking about Poisson structures we may (and will) even
identify $\CC(A)$ with the full symmetric algebra $\oplus_{n=0}^\infty
\Symm\bigl((A^*)^{\otimes n}\bigr)$. The algebra $\CC(A \times A)$ is
then identified with $\CC(A) \otimes \CC(A)$.
\begin{Definition}
We say that Poisson bracket $\delta^* : \CC(A) \wedge \CC(A) \to \CC(A)$ is
a bracket of {\em degree $N$} if
%*
\begin{displaymath}
\delta^*(A^* \wedge A^*) \subset \Symm\bigl((A^*)^{\otimes N}\bigr).
\end{displaymath}
%*
\end{Definition}
In other words, in natural local coordinates $x^i$ on the space $A$ the
Poisson tensor $\pi^{ij}(x^1,\dots,x^n)$ of the bracket $\delta^*$ must be
a homogeneous polynomial of degree $N$. It is easy to see that the values
of bracket $\delta^*$ on the space $A^* \wedge A^*$ (i.e. on linear
functions) uniquely determine, by Leibnitz identity, its values on the
whole $\CC(A) \wedge \CC(A)$, so that we will denote a restriction of
$\delta^*$ to $A^* \wedge A^*$ by the same symbol. We chose an unusual name
$\delta^*$ for the bracket because we will soon make an active use of its
dual.

Brackets of degree $1$ will be called linear, of degree $2$, quadratic,
etc. In particular, linear (Berezin--Lie) brackets are just the brackets
compatible (in the sense of the previous Section) with addition operation
in $A$. The next-simplest case is quadratic brackets (see \cite{Skl},
\cite{BAL2}, \cite{KUP2}). They play a special role in our considerations
due to the following
\begin{lemma}\label{Quadr}
A Poisson bracket compatible with the multiplication in an algebra $A$ with
a unit is quadratic. Its Poisson tensor in the unit of $A$ is zero.
\end{lemma}
\begin{proof}
Let $e_i$ be, as usual, a basis of the algebra $A$, $u = u^ie_i$ be its
unit, and $x^i$ be a dual basis in $A^*$. Taking $f = x^i,\,g = x^j$ in the
compatibility condition (\ref{Compat}), one obtains the following identity:
%*
\begin{equation}\label{Hij}
\pi^{ij}(y*z) = a_{pq}^{i} a_{st}^{j} (z^q z^t \pi^{ps}(y) + y^p y^s
\pi^{qt}(z))
\end{equation}
%*
for any $y = y^ke_k$ and $z = z^ke_k$. We are now to extract information
from (\ref{Hij}) assigning special values to $y$ and $z$.

First, take $y = z = tu$ where $t \in \Real$ is arbitrary. This gives:
%*
\begin{eqnarray*}
\pi^{ij}(t^2 u) &=&  t^2 a_{pq}^{i} a_{st}^{j} (u^q u^t \pi^{ps}(tu) + u^p
u^s \pi^{qt}(tu)) \\
&=& t^2 \bigl\langle x^i,e_p u \bigr\rangle \bigl\langle x^j, e_s
u\bigr\rangle \pi^{ps}(tu) + t^2 \bigl\langle x^i, u e_q \bigr\rangle
\bigl\langle x^j, u e_t\bigr\rangle u^p u^s \pi^{qt}(tu) \\
&=& t^2 \delta_p^i \delta_s^j \pi^{ps}(tu) + t^2 \delta_q^i \delta_t^j
\pi^{qt}(tu) \\
&=& 2t^2 \pi^{ij}(tu).
\end{eqnarray*}
%*
Denote $q(t) \bydef \pi^{ij}(tu)$. Thus, the function $q$ satisfies the
equation
%*
\begin{equation}\label{QOfT}
q(t^2) = 2 t^2 q(t).
\end{equation}
%*
Substitution $t = 0,\,1$ gives:
%*
\begin{equation}\label{QOf01}
q(0) = q(1) = 0,
\end{equation}
%*
proving the second assertion of the lemma.
Take now some $\alpha \in \Real,\ \alpha > 0$. The sequence $\alpha_n
\bydef \alpha^{1/2^n}$ tends to $1$ as $n \to \infty$, so, by continuity,
$q(\alpha_n)$ should also tend to $0 = q(1)$ (recall that, by definition,
the Poisson tensor $\pi^{ij}(x) = \{x^i,x^j\}$ is supposed to be a smooth
function). But (\ref{QOfT}) implies, by induction in $n$, that
%*
\begin{equation}\label{QOfAN}
q(\alpha_n) = \frac{q(\alpha)}{2^n \alpha^{2 - 1/2^n}}.
\end{equation}
%*
Now, combining (\ref{QOfAN}) with (\ref{QOf01}) and the fact that $q(s)$ is
smooth, one obtains the following equation:
%*
\begin{equation}\label{DerivAt1}
q'(1) = \lim_{n \to \infty} \frac{q(\alpha_n)}{1 - \alpha_n} =
\frac{q(\alpha)}{\alpha^2 \log{\alpha}}
\end{equation}
%*
and, therefore, it should be $q(\alpha) = C \alpha^2 \log{\alpha}$ for some
constant $C$ and for all $\alpha > 0$. But $q(t)$ is smooth at $t = 0$, so
that $C = 0$, and $q(t) \equiv 0$ for all $t \geq 0$. But equation
(\ref{QOfT}) shows that function $q$ is even: $q(-t) = q(t)$, and therefore
allows to extend the latter identity to all real $t$:
%*
\begin{equation}\label{HOfU}
\pi^{ij}(tu) = q(t) \equiv 0.
\end{equation}
%*

Take now in (\ref{Hij}) $z = tu$ and $y \in A$ arbitrary to obtain
%*
\begin{eqnarray*}
\pi^{ij}(ty) &=& t^2 a_{pq}^{i} a_{st}^{j} (u^q u^t \pi^{ps}(y) + y^p y^s
\pi^{qt}(tu)) \\
&=& t^2 \bigl\langle x^i, e_p u\bigr\rangle \bigl\langle x^j, e_s u
\bigr\rangle \pi^{ps}(y) \\
&=& t^2 \delta_p^i \delta_s^j \pi^{ps}(y) = t^2 \pi^{ij}(y).
\end{eqnarray*}
%*
A smooth function satisfying this identity is necessarily a quadratic form.
\end{proof}

Notice a key role that {\em smoothness} of $\pi^{ij}(x)$ plays in the
proof; see Example \ref{Iso} below.

Thus, being interested in brackets compatible with algebras, we may
restrict ourselves to quadratic ones. A general coordinate expression for
the quadratic Poisson (and also pre-Poisson) bracket $\delta^*$ is:
%*
\begin{equation}\label{QuadCoord}
\{x^i,x^j\} = c_{kl}^{ij}x^kx^l,
\end{equation}
%*
where symbol $c_{kl}^{ij}$ is skew-symmetric with respect to the upper
indices, $c_{kl}^{ij} = -c_{kl}^{ji}$. As to the lower indices, the
question of symmetry is a priori meaningless because only the sums
$c_{kl}^{ij} + c_{lk}^{ij}$ are fixed, and values of individual
$c_{kl}^{ij}$ may be chosen. A dual mapping $\delta :  \Symm(A \otimes A)
\to A \wedge A$ is given by
as
%*
\begin{equation}\label{DualCoord}
\delta(e_k \otimes e_l + e_l \otimes e_k) = (c_{kl}^{ij} + c_{lk}^{ij}) e_i
\otimes e_j.
\end{equation}
%*
(linear functions $x^i \in A^*$ are supposed to be the basis dual to $e_i
\in A$).

A tensor square $A \otimes A$ of the algebra $A$ can also be given an
algebra structure by a componentwise multiplication. Then the set of
symmetric tensors $\Symm(A \otimes A) \subset A \otimes A$ is its
subalgebra, while the set of skew-symmetric tensors $A \wedge A \subset A
\otimes A$ is a $\Symm(A \otimes A)$-bimodule. A linear mapping $D: {\frak
A} \to V$ from algebra $\frak A$ to a $\frak A$-bimodule $V$ is called a
{\em differentiation} if it satisfies the condition
%*
\begin{equation}\label{DefDif}
D(p \cdot q) = p D(q) + D(p) q
\end{equation}
%*
for all $p,\,q \in \frak A$.

\begin{theorem}\label{Diff}
Quadratic Poisson bracket $\delta^* \colon A^* \wedge A^* \to \Symm(A^*
\otimes A^*)$ is compatible with the algebra structure in $A$ if and only
if its dual mapping $\delta$ is a differentiation of $\Symm(A \otimes A)$
with values in $A \wedge A$.
\end{theorem}

The proof is just a computation: one has to check that condition
(\ref{DefDif}) for the mapping (\ref{DualCoord}) and condition of
compatibility of the mapping (\ref{QuadCoord}) with multiplication
(\ref{StructConst}) give rise to the same identities. The details are
carried out in \cite{BAL2}. The proof does not use Jacobi identity and
therefore applies for general pre-Poisson brackets as well.

Let $V$ be a linear space and $P$ be a linear operator acting from $V
\otimes V$ to itself.  Then $P^{12}$ will mean an operator $V \otimes V
\otimes V \to V \otimes V \otimes V$ acting as $P$ on the first and second
tensor component, and as the identity operator on the third one. Notations
$P^{23}$ and $P^{13}$ have similar meaning. {\em Schouten bracket}
$\SchB{P}$ is then defined as
%*
\begin{equation}\label{DefSchou}
[P^{12},P^{13}] + [P^{12},P^{23}] + [P^{13},P^{23}].
\end{equation}
%*
with $[\cdot,\cdot]$ meaning an ordinary commutator of operators.

Let $\delta^* \colon A^* \wedge A^* \to \Symm(A^* \otimes A^*)$ be a
quadratic pre-Poisson bracket, $\delta : \Symm(A \otimes A) \colon \to A
\wedge A$ be its dual, and $\widetilde \delta:  A \otimes A \to A \otimes
A$ be an arbitrary linear extension of $\delta$. Choosing this extension
means simply that we fix coefficients $c_{kl}^{ij}$ in (\ref{QuadCoord})
while a priori only sums $c_{kl}^{ij} + c_{lk}^{ij}$ are fixed.
\begin{theorem}\label{Jacobi}
Bracket $\delta^*$ is Poisson (i.e. satisfies Jacobi identity)
if and only if
%*
\begin{equation}\label{Schou}
\SchB{\widetilde \delta}(X) = 0
\end{equation}
%*
for any fully symmetric tensor $X \in A \otimes A \otimes A$.
\end{theorem}
\begin{Remark}
In particular, it is asserted in Theorem that $\SchB{\widetilde \delta}(X)$
for fully symmetric $X$ depends on $\delta$ only, and not on its specific
extension $\widetilde\delta$.
\end{Remark}

The proof is, again, a straightforward (though rather tedious) computation:
one has to check that equation (\ref{Schou}) applied to an arbitrary fully
symmetric tensor $X = X^{pqr}e_p \otimes e_q \otimes e_r$ means just the
same as Jacobi identity for the bracket (\ref{QuadCoord}) applied to an
arbitrary triple of linear functions $x^s,\,x^u,\,x^v$. Then, observe that
due to Leibnitz identity the Jacobi identity is satisfied for arbitrary
functions if and only if it is true for linear ones. By the way, the
computation does not require $A$ to have a unit.

Note also that it is not assumed in Theorem \ref{Jacobi} that bracket
$\delta^*$ is compatible with multiplication in $A$, that is, $\delta$ is a
differentiation. Assume however that $\delta$ is a differentiation, and
even has a form:
%*
\begin{equation}\label{IntDiff}
\delta(x) = \ad_r(x) \bydef [r,x]
\end{equation}
%*
for some $r = r^{ij} e_i \otimes e_j \in A \wedge A$ (it is easy to see
that all the mappings (\ref{IntDiff}) are differentiations $\Symm(A \otimes
A) \to A \wedge A$). Suppose that $A$ has a unit $u$, then the following
notation is standard: $r^{12},\, r^{23}$, and $r^{13}$ mean elements of the
tensor cube $A \otimes A \otimes A$ equal to $r^{ij} e_i \otimes e_j
\otimes u,\,r^{ij} u \otimes e_i \otimes e_j$, and $r^{ij} e_i \otimes u
\otimes e_j$, respectively. Schouten bracket $\SchB{r}$ is still defined by
expression (\ref{DefSchou}) with $[\cdot,\cdot]$ meaning now a commutator
of elements of $A \otimes A \otimes A$ rather than operators. Substitution
of (\ref{IntDiff}) into (\ref{Schou}) leads to the following theorem:

\begin{theorem}\label{MYBE}
Differentiation $\delta$ satisfies condition (\ref{Schou}) (and therefore
defines a Poisson bracket) if and only if $\SchB{r}$ is $\ad$-invariant,
i.e. commutes with any tensor of the type $S_a = u \otimes u \otimes a + u
\otimes a \otimes u + a \otimes u \otimes u,\ a \in A$.
\end{theorem}
\begin{proof}
It is easily checked that $\ad_{[a,b]} = [\ad_a,\ad_b]$ and therefore
$\SchB{\ad_r} = \ad_{\SchB{r}}$. Thus Theorem \ref{Jacobi} means that, as
soon as $\delta = \ad_r$ satisfies Jacobi identity, $\ad_{\SchB{r}}(X) = 0$
for any fully symmetric $X \in A \otimes A \otimes A$, and thus, for any
tensor of the type $S_a,\ a \in A$. Prove now the converse. If an element
$\SchB{r}$ commutes with two elements $x$ and $y \in A \otimes A \otimes
A$, then it also commutes with their sum and their product. So, if it
commutes with all the elements $S_x$, it also commutes with any element
%*
\begin{eqnarray*}
a \otimes b \otimes c + b \otimes c \otimes a + c \otimes a \otimes b + a
\otimes c \otimes b + c \otimes b \otimes a + a \otimes c \otimes b \\
= S_a S_b S_c - S_{ab} S_c - S_{ac} S_b - S_{bc} S_a + S_{acb} + S_{bca},
\end{eqnarray*}
%*
and therefore, by linearity, with any fully symmetric tensor $X \in A
\otimes A \otimes A$.
\end{proof}
\begin{corollary}\label{SubAlg}
Let $A$ be associative algebra, and $\frak G$ be a Lie subalgebra of $A_L$.
Then any element $r \in {\frak G} \wedge {\frak G}$ satisfying classical
Yang--Baxter equation (\ref{CYBE}) defines a quadratic Poisson bracket
compatible with $A$.
\end{corollary}

\section{Quadratic brackets and Poisson Lie groups}\label{Groups}

\begin{Definition}
{\em Lie bialgebra} is a Lie algebra $\frak G$ together with a linear
mapping $\Delta : {\frak G} \to {\frak G} \wedge {\frak G}$ such that:
\begin{enumerate}
\item A dual mapping $\Delta^* : {\frak G}^* \wedge {\frak G}^* \to {\frak
G}^*$ equips ${\frak G}^*$ with a Lie algebra structure.

\item The following compatibility condition is satisfied:
%*
\begin{equation}\label{Cocycle}
\Delta([a,b]) = \ad_a \Delta (b) - \ad_b \Delta (a)
\end{equation}
%*
for arbitrary $a,\,b \in {\frak G}$. Here $\ad$ stands for the adjoint
representation of the Lie algebra $ {\frak G}$ in $   {\frak G} \wedge
{\frak G}$.
\end{enumerate}
\end{Definition}

Drinfeld in \cite{DR1} proved the following theorem:
\begin{theorem}\label{Cobound}
Mapping $\Delta_r(x) = [r,x \otimes 1 + 1 \otimes x]$ where $r$ is some
element of $\frak G \wedge G$, $x \in \frak G$ and $1$ is a formal unit,
always satisfies (\ref{Cocycle}). This mapping defines a Lie bialgebra
structure if and only if $\SchB{r}$ is $\ad$-invariant.
\end{theorem}
The Lie bialgebras with such $\Delta$ are called {\em coboundary}.

Let now $G$ be a Lie group equipped with a Poisson bracket, and $\frak G$ a
corresponding Lie algebra identified with the tangent space to the group
$G$ in its unit $e$. Take $a,\,b \in {\frak G}^*$, and choose arbitrary
functions $f(x),\,g(x)$ on $G$ such that $df(e) = a,\,dg(e) = b$. Define
now the mapping $\Delta^* : {\frak G}^* \wedge {\frak G}^* \to {\frak G}^*$
by the formula
%*
\begin{equation}\label{BiAlg}
\Delta^*(a \wedge b) = d\{f,g\}(e).
\end{equation}
%*

The following theorem is due to Drinfeld (see \cite{DR1}):
\begin{theorem}\label{DrinBiAlg}
The previous definition is sound, i.e. $\Delta^*(a \wedge b)$ depends on
$a$ and $b$ only and not on a specific choice of functions $f$ and $g$. The
bracket $\{\cdot,\cdot\}$ is compatible with the multiplication in $G$
(i.e. $G$ is a Poisson Lie group) if and only if $({\frak G}, \Delta)$ is a
Lie bialgebra. This Lie bialgebra determines the corresponding Poisson
bracket uniquely.
\end{theorem}

Consider now an associative algebra $A$ with a unit $u$, and let $G$ be a
group of its invertible elements. $G$ is an open submanifold of $A$, and
therefore, if $A$ bears a Poisson bracket, then $G$ is a Poisson
submanifold (i.e., bracket can be restricted to $G$). Restriction of
quadratic bracket considered in the previous Section, gives
\begin{theorem}\label{PoiLie}
Let $A$ be an associative algebra with the unit $u$, and $\delta^* : A^*
\wedge A^* \to \Symm(A^* \otimes A^*)$ be a quadratic Poisson bracket
compatible with it. The restriction of $\delta^*$ to the group $G$ of
invertible elements of the algebra equips $G$ with a Poisson Lie structure.
The corresponding Lie bialgebra structure is given by the formula
%*
\begin{equation}\label{CoMult}
\Delta(x) = 2\delta(x \otimes u + u \otimes x).
\end{equation}
%*
\end{theorem}
\begin{proof}
Properties 1--3 from definition of Poisson bracket are obviously preserved
because $G$ is an open submanifold of $A$. Since $G$ inherits its
multiplication from $A$, the bracket and the multiplication are still
compatible, so that $G$ is a Poisson Lie group. To obtain bialgebra
structure from Drinfeld's construction, let us take functions $f$ and $g$
linear. Thus, one has really to compute a differential of the Poisson
tensor (\ref{QuadCoord}) in the point $x = u$, which obviously gives
(\ref{CoMult}).
\end{proof}

So, $\Delta^*$ equips $A^*$ with the structure of Lie algebra. Coordinate
expression for Lie bracket in the basis $x^i \in A^*$ dual to the basis
$e_i \in A$ looks like:
%*
\begin{equation}\label{DualLie}
\Delta^*(x^i \wedge x^j) = (c_{kl}^{ij} + c_{lk}^{ij})u^kx^l
\end{equation}
%*
where bracket $\delta^*$ is given by (\ref{QuadCoord}), and $u = u^ke_k$.

This formula has an interesting byproduct. Consider a mapping $\Delta_a : A
\to A \wedge A$ acting as
%*
\begin{equation}\label{LinBrack}
\Delta_a(x) = \delta(x \otimes a + a \otimes x).
\end{equation}
%*
where $a$ is an arbitrary element of $A$. It turns out to be that
$\Delta_a^* : A^* \wedge A^* \to A^*$ is always a Lie bracket, i.e.
satisfies Jacobi identity. Indeed, take $a = a^ie_i$, and write coordinate
expressions for the operation $\Delta_a^*$. One can easily see that it is
just (\ref{DualLie}) with $a_i$ substituted for $u_i$. We know from Theorem
\ref{PoiLie} that $\Delta_u^*$ satisfies Jacobi identity. But if we checked
it in coordinates we would not make use of the fact that $u$ is a unit of
algebra, and therefore the proof applies to an arbitrary $a$ as well. Note
also that formula (\ref{LinBrack}) may be viewed as an expression for the
Poisson tensor of a {\em linear} pre-Poisson bracket, and the last
assertion means then that this bracket is really Poisson. So, any quadratic
bracket compatible with an associative algebra $A$ equips it also with a
series of linear Poisson brackets.

Moreover, the following theorem holds:
\begin{theorem}\label{CompBrack}
Brackets $\Delta^*$ and $\Delta_u^*$ (where $u$ is a unit of algebra $A$)
are compatible with one another, i.e. their arbitrary linear combination
$\alpha \Delta^* + \beta \Delta_u^*$ is also a Poisson bracket.
\end{theorem}
\begin{proof}
The case $\alpha = 0$ has already been considered, so take $\alpha = 1$ and
$\beta = t$. If the function $\pi^{ij}(x) = c_{kl}^{ij}x^kx^l$ satisfies
(\ref{JaPoi}) (which means that Jacobi identity holds), then so does the
function
%*
\begin{displaymath}
\Pi^{ij}(x) \bydef \pi^{ij}(x + tu) = c_{kl}^{ij} x^kx^l + t(c_{kl}^{ij} +
c_{lk}^{ij})u^i x^j + t^2c_{kl}^{ij} u^ku^l
\end{displaymath}
%*
The last term vanishes by Lemma \ref{Quadr}, and the rest is exactly the
Poisson tensor for $\Delta^* + t\Delta_u^*$.
\end{proof}

Take up again to coboundary Lie bialgebras to prove

\begin{theorem}\label{Main}
Let $G$ be a connected simply connected Poisson Lie group such that its Lie
algebra ${\frak G} = A_L$ where $A$ is an associative algebra with the
unit. Let the corresponding Lie bialgebra be coboundary. Then $G$ contains
a discrete subgroup $\Gamma$ such that the factorgroup $G/\Gamma$ is a
Poisson Lie group, natural projection $p : G \to G/\Gamma$ is a Poisson
mapping, and $G/\Gamma$ bears a global coordinate system in which its
Poisson tensor is quadratic.
\end{theorem}
\begin{proof}
Let $({\frak G}, \Delta_r)$, where $\Delta_r$ is as in Theorem
\ref{Cobound}, be a coboundary Lie bialgebra corresponding to the Poisson
Lie structure on $G$. Then by Theorem \ref{DrinBiAlg} and Theorem
\ref{MYBE} $\ad_r^*/2$ is a quadratic Poisson bracket compatible with the
associative algebra $A$. Let $A_{inv}$ be a connected component of the unit
in the group of all invertible elements of $A$. Then the bracket
$\ad_r^*/2$ defines (by Theorem \ref{PoiLie}) a Poisson Lie structure on
the group $A_{inv}$ quadratic in the natural linear coordinates on this
group.

A Lie algebra corresponding to $A_{inv}$ is $A_L = \frak G$. Since $G$ is
connected and simply connected, then $A_{inv} = G/\Gamma$ for some discrete
subgroup $\Gamma \subset G$, and the canonical homomorphism $p : G \to
G/\Gamma$ is a covering. It allows us to define a Poisson bracket
$\{\cdot,\cdot\}$ on $G$ as an inverse image of the Poisson Lie bracket on
$A_{inv}$. We will see now that this bracket is exactly the original one,
which will complete the proof.

Define the mapping $\Delta' : \frak G \to G \wedge G$ by Drinfeld's recipe
(see (\ref{BiAlg}) above) using the bracket $\{\cdot,\cdot\}$.  Since $p$
is a covering, sufficiently small neighborhoods of units of the group $G$
and of the group $G/\Gamma = A_{inv}$ are equivalent as Poisson manifolds.
Thus $\Delta'$ should be the same mapping as the above construction would
give if applied to the group $A_{inv}$. Theorem \ref{PoiLie} gives now that
$\Delta' = \Delta_r$.

In particular, $({\frak G}, \Delta') = ({\frak G}, \Delta_r)$ is a Lie
bialgebra and therefore by Theorem \ref{DrinBiAlg} $(G,\{\cdot,\cdot\})$ is
a Poisson Lie group. Moreover, this Lie bialgebra coincides with the Lie
bialgebra corresponding to the original Poisson Lie structure on $G$, and
thus Theorem \ref{DrinBiAlg} shows that these two structures are the same.
\end{proof}

Theorem \ref{Main} may fail to be true for a non-coboundary brackets ---
consider a linear bracket on the (additive) Lie group $\Real^n$. Note also
that the proof in fact gives more than the formulation promises: quadratic
bracket on the group $G/\Gamma$ is defined by an explicit (and easy)
formula.

\section{Examples}\label{Examples}
Besides those described here, several important examples can be found in
the article \cite{STS1}.

\begin{example}\label{ab-ba}
Let $A$ be an associative algebra with a solvable $A_L$. By classical Lie
theorem any finite-dimensional linear representation of $A_L$ has an
eigenvector. Applying this to the regular representation, one obtains that
$A_L$ contains two elements, $a$ and $b$, with $[a,b] = sb$ for some
constant $s$. The element $r = a \otimes b - b \otimes a$ satisfies
Yang--Baxter equation (\ref{CYBE}) and thus, according to Corollary
\ref{SubAlg}, defines a quadratic Poisson bracket on $A$.
\end{example}

\begin{example}\label{Quaternions}
Consider a body $\bf H$ of quaternions as a four-dimensional
$\Real$-algebra. Then any element $r \in {\bf H} \wedge {\bf H}$ of the
type
%*
\begin{equation}\label{rquatern}
r = a\,\ii \wedge \jj + b\,\ii \wedge \kk + c\,\jj \wedge \kk
\end{equation}
%*
satisfies conditions of Theorem \ref{MYBE}. Corresponding Poisson bracket
on $\bf H$ is:
%*
\begin{eqnarray}
\{x^1,x^2\} &=& x^2(bx^3 - ax^4) + c((x^3)^2 + (x^4)^2) \label{RowFirst}\\
\{x^1,x^3\} &=& -x^3(cx^2 + ax^4) - b((x^2)^2 + (x^4)^2)\\
\{x^1,x^4\} &=& x^4(-cx^2 + bx^3) + a((x^2)^2 + (x^3)^2)\\
\{x^2,x^3\} &=& x^1(-bx^2 + cx^3)\\
\{x^2,x^4\} &=& -x^1(ax^2 + cx^4) \\
\{x^3,x^4\} &=& x^1(ax^3 - bx^4)\label{RowLast}
\end{eqnarray}
%*
where $x^1,\,x^2,\,x^3,\,x^4$ is a basis in ${\bf H}^*$ dual to
$1,\,\ii,\,\jj,\,\kk \in {\bf H}$.

A simply connected Lie group corresponding to $\bf H$ is a multiplicative
group of all nonzero quaternions. It contains a subgroup $\{z \in {\bf H},
\lnorm z\rnorm = 1\}$ isomorphic to $\mathop{\rm SU}(2)$. It can be checked
(see \cite{BAL2}) that this subgroup is a Poisson submanifold (and
therefore a Poisson Lie group) with respect to all the brackets
(\ref{RowFirst})--(\ref{RowLast}). Thus, $\mathop{\rm SU}(2)$ bears a
three-parameter family of Poisson Lie structures.
\end{example}

The next example shows what can happen if the algebra $A$ is associative
but has no unit.

\begin{example}\label{Nilpot}
Consider a Lie algebra $\frak G$ satisfying an identity $[[{\frak G},{\frak
G}], {\frak G}] = 0$ (i.e. $[[a,b],c] = 0$ for all $a,\,b,\,c \in \frak G$)
A particular case here is Heisenberg algebra, a three-dimensional Lie
algebra with generators $p,\,q$, and $z$, and relations $[p,q] = \Planck
z$, and $[z, p] = [z, q] = 0$. Then $a * b = \frac{1}{2} [a,b]$ is an
associative operation with $[a,b] = a*b - b*a$. An arbitrary element $r \in
{\frak G} \wedge {\frak G}$ commutes with all the symmetric tensors from
${\frak G} \otimes {\frak G}$ and therefore defines a zero Poisson bracket
on $\frak G$.

Consider now an algebra $A = \langle 1\rangle \oplus {\frak G}$ where $1$
is a formal unit. Then the simply connected Lie group $G$ corresponding to
$\frak G$ is an affine subspace $1 + {\frak G} \subset A$ (it is the usual
implementation of $G$ via exponents in the universal enveloping algebra of
$\frak G$; the thing is that $a^2 = 0$ for any $a \in {\frak G}$). Fix a
basis $\{e_i\} \in {\frak G}$; it defines a natural coordinate system in
the group $G$.  The above construction applied to $r = r^{ij} e_i \otimes
e_j \in {\frak G} \wedge {\frak G} \subset A \wedge A$ gives now the
following Poisson Lie bracket on $G$:
%*
\begin{equation}\label{NilpotBrack}
\{x^p,x^q\} = 2(r^{pl}a_{li}^{q} + r^{lq}a_{li}^{p})x^{i}
\end{equation}
%*
where $\{x^i\} \in {\frak G}^*$ is a basis dual to $\{e_i\}$, and
$a_{ij}^k$ are structure constants of the associative operation $*$ on
$\frak G$. Note that the bracket is nonzero (unlike bracket on $\frak G$
itself) and linear.
\end{example}

At last, give an example illustrating how important in Lemma \ref{Quadr} is
the fact that $\pi^{ij}(x)$ is smooth.
\begin{example}\label{Iso}
Consider the set $\Real^2$ with componentwise addition and multiplication.
It is a two-dimensional $\Real$-algebra with the unit $u = (1,1)$. Then the
Poisson bracket (singular)
%*
\begin{equation}\label{Log}
\{f(x,y),g(x,y)\} = xy\log\lmod x\rmod (\partial f/\partial x \cdot
\partial g/\partial y - \partial f/\partial y \cdot \partial g/\partial x)
\end{equation}
%*
is compatible with the multiplication in the algebra. The group of
invertible elements here is $\{(x,y) \mid x \ne 0 \ne y\}$. Its connected
component of a unit is $\Real_+^2$ where $\Real_+$ is a multiplicative
group of positive real numbers. It is isomorphic to the {\em additive}
group $\Real^2 = \{(\xi,\eta) \mid \xi,\,\eta \in \Real\}$. The
isomorphism is given by
%*
\begin{eqnarray*}
\xi = \exp(x), \quad \eta = \exp(y)
\end{eqnarray*}
%*
and maps bracket (\ref{Log}) to the linear Poisson bracket
%*
\begin{equation}\label{Lin}
\{f(\xi,\eta),g(\xi,\eta)\} = \xi(\partial f/\partial \xi \cdot \partial
g/\partial \eta - \partial f/\partial \eta \cdot \partial g/\partial \xi).
\end{equation}
%*
\end{example}

\vspace{2cm}
\noindent{\normalsize School of Mathematics and Statistics}\\
\noindent{\normalsize University of Sydney} \\
\noindent{\normalsize Sydney , AUSTRALIA, 2006}

\vspace{2cm}
\noindent{\normalsize Department of Mathematics}\\
\noindent{\normalsize Technion --- Israel Institute of Technology} \\
\noindent{\normalsize 32000, Haifa, ISRAEL}


\begin{thebibliography}{10}

\bibitem{DR1}
V.G. Drinfeld.
\newblock Hamiltonian structures on {L}ie groups, {L}ie bialgebras and
  geometric meaning of classical {Y}ang--{B}axter equations.
\newblock {\em Dokl. Akad. Nauk SSSR}, 268:285--287, 1983.

\bibitem{BAL1}
A.A. Balinsky.
\newblock Generalized quantization scheme for central extensions of lie
  algebras.
\newblock {\em J. Phys. A: Math. Gen.}, 27:L87--L90, 1994.

\bibitem{KUP3}
B.A. Kupershmidt.
\newblock A $q$-analogue of the dual space of the {L}ie algebras
  $\mathop{gl}(2)$ and $\mathop{sl}(2)$.
\newblock {\em J. Phys. A: Math. Gen.}, 26:L1--L4, 1993.

\bibitem{KUP4}
B.A. Kupershmidt.
\newblock All quantum group structures on the supergroup $\mathop{GL}(1 | 1)$.
\newblock {\em J. Phys. A: Math. Gen.}, 26:L251--L256, 1993.

\bibitem{KUP5}
B.A. Kupershmidt.
\newblock Poisson relations between minors and their consequences.
\newblock Preprint, 1994.

\bibitem{STS1}
M.A. Semenov-Tyan-Shanskii.
\newblock What is classical {$R$}-matrix ?
\newblock {\em Funct. Anal. Appl.}, 17(17):259--272, 1983.

\bibitem{Manin}
Yu. Manin.
\newblock Quantum groups and non-commutative geometry.
\newblock Montreal preprint, 1988.

\bibitem{DR}
V.G. Drinfeld.
\newblock Quantum groups.
\newblock In {\em Proc. ICM86}, volume~1, pages 798--820, Berkeley, 1987.

\bibitem{Skl}
E.K. Sklyanin.
\newblock Some algebraic structure connected with the {Y}ang--{B}axter
  equation.
\newblock {\em Funct. Anal. Appl.}, 17:263--270, 1982.

\bibitem{BAL2}
A.A. Balinsky and Yu.M. Burman.
\newblock Quadratic poisson brackets compatible with an algebra structure.
\newblock {\em J. Phys. A: Math. Gen.}, 27:L693--L696, 1994.

\bibitem{KUP2}
B.A. Kupershmidt.
\newblock Conformally multiplicative poisson structures in linear algebra and
  algebraic geometry.
\newblock In J.~Harnad and J.E. Marsden, editors, {\em Hamiltonian Systems,
  Transformation Groups and Spectral Transform Methods}, pages 177--188, 1990.

\end{thebibliography}
\end{document}